\documentclass [a4paper,notitlepage,11pt]{article}
\usepackage[pdftex]{graphicx}
\usepackage{graphicx}
\usepackage{float}
\usepackage{color}
\usepackage{hyperref}
\hypersetup{colorlinks=true,
            linkcolor={blue},
            citecolor={green}}
\usepackage[shortlabels]{enumitem}

\usepackage{amsmath, amsfonts, latexsym}
\usepackage{mathabx}
\usepackage{mathrsfs,amssymb}
\usepackage{version}
\usepackage{bm}
\usepackage{bbm}
\usepackage{flafter}
\usepackage{ulem}
\usepackage{dsfont}
\usepackage{amsthm}
\usepackage[english]{babel}
\numberwithin{equation}{section}

\newcommand{\cD}{\mathcal{D}}

\newcommand{\cF}{\mathcal{F}}

\newcommand{\cK}{\mathcal{K}}

\newcommand{\cR}{\mathcal{R}}
\newcommand{\cS}{\mathcal{S}}
\newcommand{\cT}{\mathcal{T}}

\newcommand{\cX}{\mathcal{X}}

\newcommand{\E}{\mathbb{E}}
\newcommand{\F}{\mathbb{F}}

\renewcommand{\P}{\mathbb{P}}

\newcommand{\Q}{\mathbb{Q}}
\newcommand{\R}{\mathbb{R}}

\def \proof{{\noindent \bf Proof. }}
\def \eproof{\hbox{ }\hfill$\Box$}

\newcommand{\ud}{\mathrm{d}}

\newcommand{\1}{{\bf 1}}
    
\newcommand{\set}[1]
    {\ensuremath{\{ #1 \}}}
    
\newcommand{\HP}[1] 
    {\ensuremath{\mathscr{H}^{#1}}}

\newcommand{\p}{\mathbb P}

\newcommand{\M}{\mathbb M}
\newcommand{\s}{\mathbb S}
\newcommand{\defeq}{\mathrel{\mathop:}=}

\newcommand{\diag}[1]{{\rm diag}\left[#1\right]}

\usepackage{bbm}

\DeclareMathOperator{\interior}{int}

\DeclareMathOperator{\cl}{cl}
\DeclareMathOperator{\Tr}{Tr}

\DeclareMathOperator{\essinf}{essinf}

\DeclareMathOperator{\esssup}{esssup}
\newcommand{\Rpp}{{\mathcal O}_{\!+}^{d}}

 \theoremstyle{plain}
 \newtheorem{theo}{Theorem}[section]
 \newtheorem{lem}{Lemma}[section]
 \newtheorem{prop}{Proposition}[section]
 \newtheorem{cor}{Corollary}[section]

 \newtheorem{rem}{Remark}[section]

 \theoremstyle{definition}

 \def \proof{{\noindent \bf Proof. }}
 \def \eproof{\hbox{ }\hfill$\Box$}
 \def\Im{{\rm Im}}

 \theoremstyle{remark}
 
 \newtheorem{assump}{Assumption}[subsection]

 \newcommand{\bT}
    {\ensuremath{ \mathbb{T}}}

 \newcommand{\lf}[1]{ #1^\sharp} 

 \newcommand{\pr}[1]{ #1^\flat} 

 \newcommand{\SmD}{\textcolor{black}{\cS\setminus\cD}}

\newcommand{\mage}[1]{\textcolor{black}{#1}}

\newcommand{\yellow}[1]{}
\newcommand{\bleu}[1]{\textcolor{black}{#1}}

 \title{\LARGE A comparison principle for PDEs arising in approximate hedging problems: application to Bermudan options.\footnote{The final publication is available at Springer via http://dx.doi.org/[$10.1007/s00245-017-9413-5$].}
 \footnote{This research has been sponsored by the Natixis Foundation for Quantitative Finance.}}
 \author{
         G\'{e}raldine Bouveret
         \\\small Department of Mathematics
         \\\small Imperial College London
         \\\small \sf g.bouveret11@imperial.ac.uk
         \and
         Jean-Fran\c{c}ois Chassagneux
         \\\small Laboratoire  de Probabilit\'{e}s  et Mod\`{e}les  Al\'{e}atoires
         \\\small CNRS, UMR 7599,  Universit\'{e} Paris Diderot
         \\\small \sf jean-francois.chassagneux@univ-paris-diderot.fr
}

\begin{document}

 \maketitle

 \begin{abstract}
In a Markovian framework, we consider the problem of finding the minimal initial value of a controlled process allowing to reach a stochastic target with a given level of expected loss. This question arises typically in approximate hedging problems.
The solution to this problem has been characterised by Bouchard, Elie and Touzi in \cite{BET09} and is known to solve an Hamilton-Jacobi-Bellman PDE with discontinuous operator.
In this paper, we prove a comparison {theorem} for the corresponding PDE by showing first that it can be rewritten using a continuous operator, in some cases. As an application, we then study the quantile hedging price of Bermudan options in the non-linear case, pursuing the study initiated in \cite{BBC14}.
%
    \end{abstract}

 \textbf{Keywords:} stochastic target problems, comparison principle, 
 quantile hedging, Bermudan options. 

 \vspace{1em}

 {\noindent {\sl AMS 2010 Subject Classification:} Primary, 49L25, 60J60, 93E20 ; Secondary, 49L20,  35K55.}

  \section{Introduction}
\textcolor{white}{aa}
   On a filtered probability space, $(\Omega,(\cF_t)_{t \ge 0},\P)$, and given a terminal date $T>0$, we consider two diffusion processes $\{X^{t,x}_s,\,t\leq s\leq T\}$ and $\{Y^{t,x,y,\nu}_s,\,t\leq s\leq T\}$  with values in respectively $(0,\infty)^d$ and $\R$ and initial conditions $(t,x) \in [0,T]\times (0,\infty)^d$ and $(t,x,y) \in [0,T]\times(0,\infty)^d \times \R$. The process $Y^{t,x,y,\nu}$ is controlled by $\nu$, a square integrable and progressively measurable process valued in $\R^d$.
     We are interested in the problem of finding the minimal initial value of a controlled process that allows to reach a  target with a given level of expected loss, i.e.
     \begin{equation}\label{eq Pbintro}
      v(t,x,p) \defeq \inf \left\{\begin{matrix}y\in\R^+\,\text{s.t.}\,\E\left[\ell\circ G\left(X^{t,x}_T,Y^{t,x,y,\nu}_T\right)\right]\geq p,\\\text{for some admissible controls}\,\nu\in\R^d\end{matrix}\right\}\,,
    \end{equation}
    with \textcolor{black}{$p\in I\subseteq\R$}, $\ell$ a real-valued non-decreasing function and  $G$ a real-valued function such that for $x\in(0,\infty)^d$, $y\mapsto G(x,y)$ is non-decreasing and $y\mapsto \ell\circ G(x,y)$ is right-continuous. Here, $I$ is an interval given  by the closed convex hull of the image of $\ell\circ G$, namely \textcolor{black}{$I\defeq\overline{\text{conv}}\left(\ell\circ G\left((0,\infty)^d\times\R^+\right)\right)$}.


    Problem \eqref{eq Pbintro} is coined \textit{stochastic target problem with controlled loss} by Bouchard, Elie and Touzi in \cite{BET09} who considered a non-linear Markovian formulation in a Brownian diffusion setting.     Moreau in \cite{M11} and Bouchard, Elie and Reveillac in \cite{BER15}  extended their results considering respectively jump diffusions and a non-Markovian setting.
    This problem arises when optimal management decisions are based on some risk criterion given by the loss function $\ell$. The latter belongs to the class of approximate hedging problems.

    To obtain a PDE characterisation of $v$, Bouchard, Elie and Touzi in \cite{BET09} first transformed the above problem into a stochastic target one in the $\p$-a.s. sense. To do so, they introduced an additional controlled state variable $P^{\alpha}\in\textcolor{black}{I}$ a.s. coming from the use of the martingale representation theorem. This reformulation allows then to use  the geometric dynamic programming approach introduced by Soner and Touzi  \cite{ST02GEO,ST02} for a European constraint and by Bouchard and Vu in \cite{BV10} for an obstacle constraint. However the additional controlled process in the increased state is unbounded leading to singular stochastic target problems and the Hamilton-Jacobi-Bellman equation, derived from the dynamic programming principle, involves a discontinuous operator.


     The aim of this paper is to prove a comparison theorem for the PDE solved by $v$, opening the way to direct numerical methods to compute $v$. For example, one could build a convergent numerical scheme by adapting the generalised finite difference scheme defined by Bokanowski, Bruder, Maroso and Zidani in \cite{BBMZ09}. This method has to be compared to the dual algorithm proposed in \cite{BBC14}, involving the computation of Fenchel transforms.
     {We are indeed able to prove the comparison theorem under a setting involving a semi-linear dynamics for $Y$ and for unconstrained controls $\nu$}. To the best of our knowledge, this is the first time that such a result is obtained in this non-linear setting. This answers a question raised in  \cite[Section 4]{bouchardbsdes} (preprint version of \cite{BER15}).
     One key step  is performed by using a renormalisation argument (see Section \ref{se altpde charac1}) to obtain a new continuous operator, in the spirit of \cite{Bru05}.
      However the new operator has a non-linearity in front of the time derivative. We therefore  rely, as e.g. in \cite{IL90}, on a strict super-solution argument to prove that a comparison result holds (see Section \ref{se strictsupersol1} and Section \ref{se CPR}).

     To the best of our knowledge no comparison theorem has been proved so far for the PDE solved by $v$ in the case where the controls $\nu$ belong to a constrained set. This case corresponds to a constraint imposed on the gradients of the solution. This interesting problem is left for further research. Let us mention however, that using some approximation argument, Bouchard and Vu \cite{BV12} are able to obtain a convergent numerical procedure to compute the value function $v$ at its continuity point, see also \cite{BD13}.

    In Section \ref{se appli}, we use our comparison result to characterise the quantile hedging price of Bermudan options in a non-linear setting, pursuing the study initiated in \cite{BBC14}. Quantile hedging problems have been introduced  by F\"{o}llmer and Leukert (see \cite{FL99}) for European-type claims and later studied by \cite{BET09, BBC14, JKT13} amongst other.
Precisely, we consider the problem of a trader who wants to find a hedging strategy $\nu$ and an initial endowment $y$ such that his hedging portfolio $Y^{t,x,y,\nu}_\cdot$ stays above a claim of the form $g(\cdot,X^{t,x}_\cdot)$ over a set of deterministic dates with a given probability. Here $X^{t,x}$ models the evolution of some risky assets, assuming that their value is $x$ at time $t$. The conservative case would be to ensure that the insurer meets a risk target over time almost surely. This is especially true in a Solvency II world where the risk should be monitored over time. However this constraint is too restrictive and we want to keep a flexible framework to allow for reasonable opportunities to make profits and this even with a limited available capital. Thus we weaken the constraint and express it in probability. 
In practice, the solvency constraint comes from an outside party, a minimal requirement for a fund manager, or the willingness to avoid a huge dis-utility. As usual in the Bermudan setting, the difficulty comes from the iteration over the time intervals of the characterisation obtained in the European case. In our case, the time-boundary condition on each interval is the most problematic issue since it involves a ``facelift'' phenomenon. Its full characterisation is obtained by using the continuity of the value function on the previous time interval, which comes from the application of the comparison {theorem}.

\vspace{1em}
The rest of the paper is organised as follows. In Section 2, we introduce our framework for the study of stochastic target problems with controlled loss. In Section 3, we prove the comparison theorem. First, we obtain a new PDE characterisation for the value function involving a continuous operator. Then, we use a strict super-solution argument to obtain the comparison result. Finally, in Section 4, we study the quantile hedging price of a Bermudan option in our non-linear setting.

    \vspace{1em}

    \textbf{Notations.} Let $d\ge 1$ be an integer. Any vector $x$ of $\R^{d}$ is seen as a column vector. We denote by $|x|$ the norm $1$ of $x$, by $\|x\|$ its norm $2$ and by $x^{\top}$ its transpose. The notation $\M^{d}$ denotes the set of $d$-dimensional square matrices and $\s^{d}$ is the subset of elements of $\M^{d}$ that are symmetric.
    We set $\mathrm M^{\top}$ the transpose of $\mathrm M\in \M^{d}$, while $\Tr[ \mathrm M]$ is its trace. We respectively denote by $\1$ and $\textbf{I}$ the $d$-dimensional unit column vector and the $d\times d$-dimensional unit  matrix. Let $\psi : (t,x,p)\in [0,T]\times\R^d\times \R \mapsto \psi(t,x,p)$. If it is smooth enough, we denote by $\partial_{t}\psi$ its derivative with respect to $t$ and by $\mathrm D\psi$ its Jacobian matrix with respect to the space variables whose rows are given by $\mathrm D_{x}\psi$ and $\mathrm D_{p}\psi$, i.e. the derivative with respect to $x$ and $p$. The Hessian matrix with respect to the space variables is $\mathrm D^{2}\psi$ whose elements are given by $\mathrm D_{xx}\psi,\mathrm D_{pp}\psi,\mathrm D_{xp}\psi,\mathrm D_{px}\psi$, i.e. the second derivative with respect to $x$ and $p$ and the cross derivatives. For a given function $f\in\R^d$, $f^{-1}$ stands for its inverse. Moreover, for a given function $f\in\R$, $\overline{\text{conv}}(f)$ is the closed convex envelope of $f$.

    \yellow{
    Furthermore we denote by $\mathcal{S}$ the sphere of $\R^{d+1}$ of radius $1$ and $\cD$ the set of vectors $\eta \in \cS$     such that their first component $\eta^1 = 0$.}

    \yellow{For a vector $\eta \in \cS \setminus \cD$, we denote by $\pr{\eta}:=\frac1{\eta^1} (\eta^2,\dots,\eta^{d+1})^\top \in \R^d$. To a vector $\alpha$ in $\R^d$ we associate a vector $\lf{\alpha}$ in $\cS \setminus \cD$, namely
    $\lf{\alpha} = (\frac1{\sqrt{1 + |\alpha|^2}}, \frac{\alpha}{\sqrt{1 + |\alpha|^2}})^\top$. We observe that $\pr{(\lf{\alpha})}= \alpha$.}

    \yellow{We will denote $\textbf{e}\defeq(0,1,...,1)^\top$ the $d+1-$dimensional column vector with all elements equal to $1$ except the first one.} Given a locally bounded map $v$ on an open subset $B$ of $\R^d$, we define the lower and upper semi-continuous envelopes
    \begin{align*}
     v_*(b)\defeq\liminf_{B\ni b'\rightarrow b} v(b')\quad&\text{and}\quad v^*(b)\defeq\limsup_{B\ni b'\rightarrow b} v(b'),\,b\in\cl(B)\,.
    \end{align*}
    Finally $\Rpp\defeq (0,\infty)^d$. All over the paper, inequalities between random variables have to be understood in the $\p$-a.s. sense.

   \section{Problem statement}\label{se pbstatement}

   In the sequel we work with a finite time horizon $T>0$. Let $\Omega$ be the space of $\R^d$-valued continuous functions $(\omega_t)_{t\leq T}$ on $[0,T]$, $d\geq 1$, endowed with the Wiener measure $\p$. We denote by $W$ the coordinate mapping, i.e. $\left(W(\omega)_t\right)_{t\leq T}$ for $\omega \in \Omega$ so that $W$ is a $d$-dimensional Brownian motion on the canonical filtered probability space $\left(\Omega,\mathcal{F},\F,\p\right)$. In the latter $\mathcal{F}$ is the Borel tribe of $\Omega$ and $\F\defeq\left\{\mathcal{F}_t,0\leq t\leq T\right\}$ is the $\p$-augmentation of the filtration generated by $W$. Let $\mathcal{U}$ be the collection of $\R^d$-valued progressively measurable processes in $L^2([0,T]\times \Omega)$. 
   For $t\in [0,T], (x,y)\in \Rpp\times\R$ and $\nu \in \mathcal{U}$ the processes $X^{t,x}$ and $Y^{t,x,y,\nu}$ are  defined as the solution to
   the following stochastic differential equations (SDE)
    \begin{align*}
      X^{t,x}_s&=x+\int_{t}^{s}\diag{X^{t,x}_r} \mu\left(X^{t,x}_r\right) \,  \ud r +\int_{t}^{s}  \diag{X^{t,x}_r}\sigma\left(X^{t,x}_r\right) \, \ud W_r \,\,,
      \\Y^{t,x,y,\nu}_s&=y+\int_{t}^{s} \mu_Y\left(r,X^{t,x}_r,Y^{t,x,y,\nu}_r,\nu_r\right) \, \ud r +\int_{t}^{s} \nu^\top_r \sigma \left(X^{t,x}_r\right) \, \ud W_r \,\,,
     \end{align*}
     where $\mu: \Rpp \rightarrow \R^d$, $\sigma: \Rpp \rightarrow \M^d$ and $\mu_Y:[0,T]\times\Rpp\times\R\times \R^d\rightarrow\mage{\R}$ are continuous functions.

     We assume that $X^{t,x}$ takes its value in $\Rpp$ when the initial condition $x$ is in $\Rpp$. \mage{Moreover the process $Y^{t,x,y,\nu}$ is valued in $\R$.}
     To guarantee that the above processes are well defined, we shall use throughout the paper the following standing assumption.
\noindent \mage{\noindent\textbf{Standing Assumption}
     \begin{enumerate}
     \item For all $(r,r',x,x',y,y',\upsilon,\upsilon')\in[0,T]^2\times(\Rpp)^2\times\R^{2}\times\R^{2d}$
      \begin{align}\label{as coeflipY}
    &|\mu_Y(r',x',y',\upsilon')-\mu_Y(r,x,y,\upsilon)|
    \nonumber\\&\quad\quad\quad\quad\quad\quad\le L\left(|r'-r|+|x'-x|(1+|\upsilon'|+|\upsilon|)+|\upsilon'-\upsilon|+|y'-y|\right),\,
    \end{align}
      for some Lipschitz constant $L>0$. 
     \item The functions $\mu_X :x\in\Rpp\mapsto\text{diag}[x]\mu\left(x\right)\in \R^d$, $\sigma_X:x\in\Rpp\mapsto
\text{diag}[x]\sigma\left(x\right)\in\s^d$ are Lipschitz continuous with
some Lipschitz constant $L>0$.
 \item The function $\sigma$ is invertible and $x\in\Rpp\mapsto\sigma^{-1} \left(x\right)\in\s^d$  is $L$-Lipschitz continuous on $\Rpp$, for some $L>0$.
    The function $\mu$, $\sigma$ and $\sigma^{-1}$ are bounded by a constant $\Lambda > 0$.
\end{enumerate}}

\mage{Finally,  we shall sometimes use the following monotonicity assumption.
\begin{assump}[Drift monotonicity]\label{as dm}
     The function $y\mapsto\mu_Y(\cdot,y,\cdot)$ is increasing.
\end{assump}
\begin{rem}
     It is well known that, in our context, setting $\tilde{Y}_t:=e^{\lambda t}Y_t$, $t \in [0,T]$, we
     obtain that the drift of $\tilde{Y}$ satisfies Assumption \ref{as dm} when \eqref{as coeflipY} is in force and provided that $\lambda$ is big enough.
\end{rem}
     To ease the notations, we will thus simply assume later in the proof of the comparison principle that Assumption \ref{as dm} is in force.
     This is made clear in Remark \ref{re ope increasing} below.
     }
%
%


Then we denote by $\mathcal{U}_{t,x,y}$ the subset of $\mathcal{U}$ for which the process  $Y^{t,x,y,\nu}_\cdot\ge 0$ on $[t,T]$.\yellow{and
\begin{gather*}\label{eq sqint}
\E\left[\int_{t}^{T} \|\nu^\top_r \sigma\left(X^{t,x}_r\right)\|^2\ud r\right]  <\infty\,,
\end{gather*} for all initial data $(t,x)\in[0,T]\times\Rpp$ and $y\ge 0$.
}

\begin{rem}\label{re exdrifty}
In a financial setting the process $X^{t,x}$ is an underlying process representing the price of some risky assets while the process $Y^{t,x,y,\nu}$ is the wealth process where the control $\nu^i$ stands {for the amount invested in assets $i$}. In our setting two typical examples of $\mu_Y$ satisfy \eqref{as coeflipY}:

(i) The usual case of linear pricing, where there is a risk premium $\zeta(x)\defeq \sigma^{-1}(x)(\mu(x)-r\mathbf{1})$ with $r$ a risk-free interest rate and  \mage{$\mu_Y(x,y,\upsilon)\defeq\zeta^\top(x)\sigma(x)\upsilon$}.

(ii) A case of non-linear pricing, coming from a market imperfection when there are two non-negative rates $R$ (the borrowing rate) and $r$ (the lending rate) with $R>r$ \cite[Example 1.1]{el1997backward}
\begin{align*}
\mu_Y(x,y,\upsilon) = ry + \zeta(x)^\top\sigma(x)\upsilon -(R-r)\left( y - \upsilon^\top \1\right)_-\,,
\end{align*}
{where $\zeta$ stands for the risk premium.}
\end{rem}

    Now let $\ell: \R\rightarrow\R$ be a non-decreasing function and $G:\R^{d+1}\mapsto\R$ be a measurable map such that for any $x\in\Rpp$, $y\mapsto G(x,y)$ is non-decreasing and $y\mapsto \ell\circ G(x,y)$ is right-continuous.
    {
    We also assume that $I\defeq\overline{\text{conv}}\left(\ell\circ G\left(\Rpp\times\R^+\right)\right)$,  the closed convex hull of the image of $\ell\circ G$, is a compact interval of $\R$. In our application below, it is clear
     that -- up to a proper rescaling of the $\ell$ function -- one can consider that $I=[0,1]$ and we shall work under this setting from
     now on.
    }
    For $(t,x,p)\in [0,T]\times\Rpp\times [0,1]$, we  then define the stochastic target problem with controlled loss as
     \begin{equation}\label{eq Pb1}
      v(t,x,p) \defeq \inf \Big\{y\in\R^+: \exists \hspace{1mm} \nu\in \mathcal{U}_{t,x,y}\hspace{2mm}\text{s.t.}\,\E\left[\ell\circ G\left(X^{t,x}_T,Y^{t,x,y,\nu}_T\right)\right]\geq p\Big\}\,.
    \end{equation}

     {\begin{assump}\label{assump polygrowth0} We assume that there exists a constant $\beta>0$  such that
     \begin{align} \label{eq polygrowth0}
     \lvert v(t,x,1)\rvert\leq \beta(1+\lvert x\rvert^{k}),\,k\ge1,\; \text{ for all } (t,x)\in[0,T]\times\Rpp\,.
     \end{align}
      \end{assump}}
 \begin{rem}\label{re polygrowth0} Since, for all $p\in[0,1]$, $0\le v(\cdot,p)\le v(\cdot,1)$ the previous assumption implies the
 condition \eqref{eq polygrowth0} holds true for $v$ on $[0,T]\times\Rpp\times[0,1]$.
 \end{rem}
     We assume that $\ell\circ G\left(X^{t,x}_T,Y^{t,x,y,\nu}_T\right)$ is square integrable for all initial conditions and for all $\nu\in\mathcal{U}_{t,x,y}$. Bouchard, Elie and Touzi proved in \cite{BET09} that, in that case, \eqref{eq Pb1} can be reduced to
       \begin{equation}\label{eq Pb1Alt}
         v(t,x,p) = \inf \Big\{y\in\R^+: \exists \hspace{1mm} (\nu,\alpha)\in \mathcal{U}_{t,x,y}\times\mathcal{A}_{t,p}\hspace{2mm}\text{s.t.}\,\ell\circ G\left(X^{t,x}_T,Y^{t,x,y,\nu}_T\right)\ge P^{t,p,\alpha}_T\Big\}\,,
       \end{equation}
       where for $p\in\bleu{[0,1]}$, $\mathcal{A}_{t,p}$ is the set of $\R^d$-valued $\F$-progressively measurable and square integrable processes $\alpha$ such that
       \begin{align*}
        P^{t,p,\alpha}_T\defeq p+\int^{T}_{t}\alpha^\top_s\ud W_s\in\bleu{[0,1]}\,.
        \end{align*}
As a consequence the problem reduction implies to work with an unbounded set of controls whatever the set of controls $\mathcal{U}$ is and then to deal with a discontinuous Hamilton-Jacobi-Bellman operator. The latter makes the proof of a comparison result difficult.
    The aim of this paper is thus to provide a comparison principle in the above framework (see Section \ref{se CP}). We will use this
    result to provide a full PDE characterisation of the quantile hedging price  of a Bermudan option in our  non-linear framework (see Section \ref{se appli}). 

      \section{A comparison principle}\label{se CP}
      In this section we prove a comparison principle for the PDE satisfied by the value function given in \eqref{eq Pb1}.
      As observed in \cite[Section 4]{bouchardbsdes} (preprint version of \cite{BER15}), this is not straightforward as this PDE is naturally obtained using a discontinuous operator, see equations \eqref{eq subsolBET}-\eqref{eq supersolBET} below.
      In a first step, we are able  to show that any solution to this PDE can be characterised by a PDE involving a continuous operator.
      However this operator is non-standard as it involves a non-linearity in the time-derivative. Nevertheless, using a strict super-solution approach, we manage to prove a comparison theorem for this new PDE (and thus the original one).

   \subsection{Alternative PDE characterisation inside the domain}\label{se altpde charac1}

   Let us start with some definitions.
 For  $ (t,x,y)\in [0,T]\times\Rpp\times\R^{+}$,
$q\defeq\left(\begin{matrix}q^{x}\\q^p\end{matrix}\right)\in\R^{d+1}$ and $A\defeq\left(\begin{matrix}A^{xx}&A^{xp}\\A^{xp^\top}&A^{pp}\end{matrix}\right)\in\s^{d+1}$, denoting $\Xi \defeq(t,x,y,q,A)$, we define
              \begin{align*}
        \mathrm F^*(\Xi)\defeq \limsup_{\Xi'\rightarrow \Xi}\mathrm F(\Xi')\,,
      \end{align*}
      with
      \bleu{
      \begin{gather*}
        \mathrm F(\Xi)\defeq\sup_{(\upsilon,a)\in\mathcal{N}(x,y,q)}
        \left\{
              \mu_Y(t,x,y,\upsilon)-\mu_X^\top(x)q^x-\frac{1}{2}\Tr\left[\bar\sigma\bar\sigma^\top(x,a)A \right]
           \right\}\,,
      \end{gather*}
      where
      \begin{align*}
      \mathcal{N}(x,y,q)\defeq \{(\upsilon,a)\in \R^d \times\R^d\,:\,\upsilon^\top \sigma(x)\defeq  q^\top\bar\sigma(x,a) \}
       \;\text{ and }\;
            \bar\sigma(x,a)\defeq \left(\begin{matrix}\sigma_X(x)\\
           a^\top \end{matrix}\right)\,,
           \end{align*}
       }
\mage{and where we recall the notations
\begin{align*}
\mu_X(x)\defeq\text{diag}[x]\mu(x)\quad\sigma_X(x)\defeq \text{diag}[x]\sigma(x)\,.
\end{align*}}
 \bleu{
 Let us observe that as $\sigma$ is invertible the previous expression can be simplified as $\upsilon$
 is then a function of the variable $a$. We thus introduce
 \begin{align*}
  J^a(\Xi) = \hat\mu_Y(t,x,y,q,a)-\mu_X^\top(x)q^x-\frac{1}{2}\Tr\left[\bar\sigma\bar\sigma^\top(x,a)A \right]\,,
 \end{align*}
 where $\hat{\mu}_Y(t,x,y,q,a)\defeq {\mu}_Y(t,x,y,   (q^\top \bar\sigma(x,a)\sigma(x)^{-1})^\top)$ and observe that
 \begin{align*}
 \mathrm F(\Xi) = \sup_{a \in \R^d} J^a(\Xi)\;.
 \end{align*}
 }


For the reader's convenience, we will write
$$\mathrm F\varphi(t,x,p) \text{ for } \mathrm F(t,x,\varphi(t,x,p),\partial_t\varphi(t,x,p),\mathrm D\varphi(t,x,p),\mathrm D^2\varphi(t,x,p))\,.$$ This writing will hold for any super-/sub-solution operator defined hereinafter.

    Bouchard, Elie and Touzi in \cite{BET09} proved that on $[0,T)\times\Rpp\times(0,1)$, $v^*$ is a viscosity sub-solution of
\begin{align}\label{eq subsolBET}
&\min\{v^*,-\partial_t\varphi+{\mathrm F}\varphi\}\le0\,,
\end{align}
and $v_*$ is a viscosity super-solution of
\begin{align}\label{eq supersolBET}
&-\partial_t\varphi+{\mathrm F}^*\varphi\ge0\,.
\end{align}
As mentioned before the problem here stems from the fact that the Hamilton-Jacobi-Bellman operator is lower semi-continuous and not upper semi-continuous.

      As a consequence we will first work towards an alternative PDE characterisation of $v$ (see Theorem \ref{th AltPDECharac0} below) that will allow us to express both the sub-solution and super-solution properties with
      \bleu{a continuous operator}.

    \mage{Now let us denote by $\mathcal{S}$ the sphere of $\R^{d+1}$ of radius $1$ and by $\cD$ the set of vectors $\eta \in \cS$  such that their first component $\eta^1 = 0$. For a vector $\eta \in \cS \setminus \cD$, we denote $\pr{\eta}:=\frac1{\eta^1} (\eta^2,\dots,\eta^{d+1})^\top \in \R^d$}. Moreover we define for $\Theta\defeq(t,x,y,b,q,A)\in\mage{[0,T]\times}\Rpp\times\R^+\times\R\times\R^{d+1}\times\s^{d+1}$
  the following operator
%
    \bleu{
        \begin{align*}
        \mathrm H^\eta(\Theta)\defeq \left \{
        \begin{array}{cl}
              (\eta^1)^2 \left(-b+ J^{\pr{\eta}}(t,x,y,q,A)\right)
              &\,, \eta \in \SmD
              \\
              -\frac{1}{2}A^{pp} &\,, \eta \in \cD
          \end{array}\,.
          \right.
      \end{align*}
    }

    Observing that for $\eta \in \SmD$ the above operator reads
     \bleu{
     \begin{gather*}
               (\eta^1)^2\left(-b+\hat{\mu}_Y(t,x,y,q,\pr{\eta})-{\mu}_X^\top(x)q^x-\frac{1}{2}\Tr\left[\sigma_X\sigma_X^\top(x)A^{xx}\right]
               -\frac{1}{2}\|\pr{\eta}\|^2A^{pp}
               -(\pr{\eta})^\top\sigma_X^\top(x)A^{xp}
               \right)
            \,,
       \end{gather*}
       }
       we can make the following remark.
\bleu{
      \begin{rem}\label{re opcont}
      It follows from \eqref{as coeflipY} that the operator $\eta \mapsto \mathrm{H}\mage{^\eta}$ is continuous
      on $\cS$, in particular,
      \begin{align*}
       \sup_{\eta \in \cS} \mathrm H^\eta(\Theta) = \sup_{\eta \in \SmD} \mathrm H^\eta(\Theta)\,.
      \end{align*}
      \end{rem}
}
      We can now state the alternative PDE characterisation of $v$.
      \begin{theo}\label{th AltPDECharac0}
          On $[0,T)\times\Rpp\times(0,1)$,
       \bleu{$v^*$ (resp. $v_*$) is a viscosity sub-solution (resp. super-solution) of
       \begin{align} \label{eq TargetSolAlt}
	    \mage{\mathcal{H}}\varphi = 0 \quad\quad \text{ with }\quad
       \mage{\mathcal{H}}(\Theta) = \min \set{\,y\,,\,\sup_{\eta\in\cS}\bleu{\mathrm{H}}^\eta(\Theta)}\,,
       \end{align}
       where $\Theta=(t,x,y,b,q,A)\in\mage{[0,T]\times}\Rpp\times\R^+\times\R\times\R^{d+1}\times\s^{d+1}$.
       }
      \end{theo}

       \proof
       \textbf{Step 1.} \textit{Proof of the sub-solution property.} \\Let $\varphi$ be a smooth function such that $\max_{[0,T)\times\Rpp\times(0,1)}(v^*-\varphi)(t,x,p)=(v^*-\varphi)(t_0,x_0,p_0)=0$.\\
       It follows from \eqref{eq subsolBET} that $$\min\{v^*(t_0,x_0,p_0),-\partial_t\varphi(t_0,x_0,p_0)+\mathrm F\varphi(t_0,x_0,p_0)\}\le 0\,.$$
       We will prove that  $$-\partial_t\varphi(t_0,x_0,p_0)+\mathrm F\varphi(t_0,x_0,p_0)\le 0\Rightarrow\sup_{\eta\in\cS}\bleu{\mathrm{H}}^\eta\varphi(t_0,x_0,p_0)\le 0\,,$$
       which will lead to the sub-solution part of \eqref{eq TargetSolAlt} as $v^*\ge 0$ by definition.

    \bleu{
     By definition of $\mathrm F$ we have that for all $a \in \R^d$,
     \begin{align*}
      - \partial_t \varphi(t_0,x_0,p_0) + J^a\varphi(t_0,x_0,p_0) \le 0\;.
     \end{align*}
     For all $\eta \in \SmD$, we then obtain
     \begin{align*}
     (\eta^1)^2 \left( - \partial_t \varphi(t_0,x_0,p_0) + J^{\pr{\eta}}\varphi(t_0,x_0,p_0) \right) \le 0\,.
     \end{align*}
   By continuity of $\eta \mapsto \mathrm H^\eta$, we thus obtain for all $\eta \in \cS$,
      \begin{align*}
       \mathrm H^{\eta} \varphi(t_0,x_0,p_0) \le 0\,.
      \end{align*}
    \mage{The arbitrariness of $\eta$} concludes the proof for this step.
    }

        \textbf{Step 2.} \textit{Proof of the super-solution property.} \\
        Let $\varphi$ be a smooth function such that
        $\min_{[0,T)\times\Rpp\times(0,1)}(v_*-\varphi)(t,x,p)=(v_*-\varphi)(t_0,x_0,p_0)=0$. We note that by definition $v_* \ge 0$ so that we just
        have to verify that $$\sup_{\eta\in\cS}\mathrm{H}^\eta\varphi(t_0,x_0,p_0)\ge 0\,.$$
       According to \eqref{eq supersolBET} we have, $$-\partial_t\varphi(t_0,x_0,p_0)+\mathrm F^*\varphi(t_0,x_0,p_0)\ge 0\,.$$
        By definition of $\mathrm F^*$ we can find sequences $t_k\in[0,T),\, (x_k,p_k)\in\Rpp\times(0,1),\,y_k\geq0,\,q_k\defeq(q_k^{x},q_k^{p})\in\R^{d+1}$, and a symmetric matrix $A_k\in\s^{d+1}$ such that
    \begin{align}
     (t_k,x_k,p_k)\rightarrow (t_0,x_0,p_0)\quad \text{and}\quad |(y_k,q_k,A_k)-(\varphi,\mathrm D\varphi,\mathrm D^2\varphi)(t_0,x_0,p_0)|\leq k^{-1}\,,\label{eq supersolinterm1}
    \end{align}
    and
    \begin{equation*}
        -{\partial_t\varphi(t_0,x_0,p_0)}+\mathrm F(t_k,{x}_k,y_k,q_k,A_k)\geq -k^{-1}\,.
    \end{equation*}
    Then we can find a maximising sequence $\mage{a_k} \in \R^d$ such that
    \begin{align*}
        -\partial_t\varphi(t_0,x_0,p_0)+J^{\mage{a_k}}(t_k,x_k,y_k,q_k,A_k) \geq -2k^{-1}\,.\label{eq supersolinterm2}
    \end{align*}

     Now consider \mage{$\eta_k := (\frac1{\sqrt{1 + \|a_k\|^2}}, \frac{a_k}{\sqrt{1 + \|a_k\|^2}})^\top$. Note that $\eta_k\in\cS \setminus \cD$.}
   Therefore we have
         \begin{align*}
        (\eta_k^1)^2\left(-\partial_t\varphi(t_0,x_0,p_0)+J^{\pr{\eta}_k}(t_k,x_k,y_k,q_k,A_k)\right) \geq -2k^{-1}(\eta_k^1)^2\,.
	\end{align*}

       Hence using the relative compactness of the set $\SmD$ we have the existence
       of a subsequence such that $\lim_{k'\rightarrow\infty}\eta_k'=\bar\eta$ with $\bar\eta\in\cS$.
       Moreover using \eqref{eq supersolinterm1}, \eqref{as coeflipY} and the standing hypotheses on the coefficients of $X$ we obtain
       $$\sup_{\eta\in\cS}\bleu{\mathrm{H}}^\eta\varphi(t_0,x_0,p_0)\geq \bleu{\mathrm{H}}^{\bar\eta}\varphi(t_0,x_0,p_0)\geq0\,,$$
       which concludes the proof for this step.
       \eproof

\begin{rem}\label{re conv} 
{\rm (1)}  We notice that 
the sub-solution property implies that $\mathrm D_{pp}\varphi(t_0,x_0,p_0)\ge 0$.\\
{\rm (2)} Note that,  if there exists a risk premium $\zeta(x)\defeq \sigma^{-1}(x)\mu(x)\in \R^d$ (i.e. Remark \ref{re exdrifty} (i) with $r=0$) such that $\mu_Y(x,y,\nu)\defeq\zeta^\top(x)\sigma(x)\nu$, \eqref{eq TargetSolAlt} implies that $v^*$ is a viscosity sub-solution of
$$\min\{v^*,\Lambda^+(\mathrm M)\}\le0\,,$$ and  that $v_*$ is a viscosity super-solution of
$$\Lambda^+(\mathrm M)\ge0\,,$$ where $\Lambda^+(\mathrm M)$ denotes the highest eigenvalue of the symmetric matrix $\mathrm M$ and where for  $\Theta\defeq(x,b,q^p,A)\in\Rpp\times\R\times\R\times\s^{d+1}$ the matrix \mage{$\mathrm M$ reads}
    \begin{gather*}
     \mathrm M(\Theta)\defeq{\mathrm M}_1(\Theta)+{\mathrm M}_2(\Theta)
     \,,
   \end{gather*}
   with ${\mathrm M}_1$ and ${\mathrm M}_2$ two $(d+1)\times(d+1)$ matrices defined as
   \begin{gather*}
     {\mathrm M}_1(\Theta)\defeq\left(
     \begin{matrix}
  -b
       -\frac{1}{2}\Tr\left[\sigma_X\sigma_X^\top(x)A^{xx}\right]& 0 & \cdots &0 \\
  0 & -\frac{1}{2}A^{pp} & \cdots & 0 \\
  \vdots  & \vdots  & \ddots & \vdots  \\
  0 & 0 & \cdots & -\frac{1}{2}A^{pp}
  \end{matrix}\right)\,,
   \end{gather*}
   and
   \begin{gather*}
     {\mathrm M}_2(\Theta)\defeq\left(
     \begin{matrix} 0&\frac{1}{2}\left[\zeta(x)q^{p}-\sigma_X(x)A^{xp}\right]^\top
     \\\frac{1}{2}\left[\zeta(x)q^{p}-\sigma_X(x)A^{xp}\right]&0
     \end{matrix}
     \right)\,.
   \end{gather*}
This is in the spirit of \cite{BBMZ09, BPZ16, Bru05}.
\end{rem}

\begin{rem} \label{re ope increasing}
By an usual change of variable argument, it is easily seen that if $v$ is a super-solution (resp. sub-solution) of
 \begin{align*}
  \sup_{\eta \in \cS} {\mathrm H}^\eta \varphi = 0\,,
 \end{align*}
then $\tilde{v}(t,x,p) := e^{\lambda t} v(t,x,p)$, for some $\lambda\mage{>0}$, is a super-solution (resp. sub-solution) of
 \begin{align*}
  \sup_{\eta \in \cS} \mathrm{\tilde{H}}^\eta \varphi = 0\,,
 \end{align*}
with
        \begin{align*}
        \mathrm{\tilde{H}}^\eta (\Theta) \defeq \left \{
        \begin{array}{cl}
              (\eta_1)^2 \left(-b+ \tilde{J}^{\pr{\eta}}(\Theta)\right)
              &\,, \eta \in \SmD
              \\
              -\frac{1}{2}A^{pp} &\,, \eta \in \cD
          \end{array}
          \right.\,,
      \end{align*}
      where
      $$
      \tilde{J}^{\pr{\eta}}(\Theta) =
      \lambda y + \mage{e^{\lambda t}}J^{\pr{\eta}}(t,x,e^{-\lambda t}y,e^{-\lambda t}q,e^{-\lambda t}A) \;.
      $$
For $\lambda > L$, we observe that $\tilde{J}^{\pr{\eta}}$ is strictly increasing in $y$. From now on, we will thus assume that $\hat{\mu}_Y$ -- and thus $J^{\pr{\eta}}$ -- is strictly increasing in $y$. Namely, we will assume that \mage{Assumption \ref{as dm}} is in force. 
\end{rem}

            \subsection{Strict super-solution property and modulus of continuity}\label{se strictsupersol1}

       The operator $\mathrm H^\eta$ has a non-linearity in front of the time-derivative. We then have to rely on a strict super-solution argument to prove that a comparison result holds for the non-linear PDE solved by $v$. This argument has been used, for example, by Ishii and Lions in \cite{IL90} and Cheridito, Soner and Touzi in \cite{CST05}. We thus have beforehand to introduce the following lemma.

\textcolor{black}{
             \begin{lem}[Strict super-solution property]\label{le StrictSupersol}
         Let us define on $[0,T]\times\Rpp\times\bleu{[0,1]}$ the smooth positive functions
          $\phi(t,p) \defeq e^{\kappa(T-t)} (\theta - \mage{\frac{e^{-4\bar{c}p}}{2}})\,$, $h(t,x) \defeq      e^{\kappa(T-t)} ( |x|^{2k}+|x|^{-2}  )$,   and
          \begin{align*}
          f(t,x,p)\defeq\phi(t,p)+h(t,x) > 0\,,
          \end{align*}
         for some $\kappa, \theta > 1$, with $\bar{c} \defeq L \times \Lambda$ \mage{(recall Standing Assumption)} and $k$ introduced in Assumption \ref{assump polygrowth0}.
\\
\vspace{5pt}
        Let $V$ be a lower semi-continuous super-solution of \eqref{eq TargetSolAlt}.
         Then, for $\kappa$ and $\theta$ big enough, the function $V+\xi f$, $\xi>0$, is a strict viscosity super-solution of \eqref{eq TargetSolAlt} on $[0,T)\times\Rpp\times(0,1)$, i.e. given a smooth function $\varphi$ such that $\min_{[0,T)\times\Rpp\times(0,1)}((V+\xi f)-\varphi)(t,x,p)=((V+\xi f)-\varphi)(t_0,x_0,p_0)=0$, one has
         \begin{align}
           \mathcal{H} \varphi(t_0,x_0,p_0)\ge \xi \varrho \,, \label{eq strict supersol}
         \end{align}
         for some $\varrho > 0$.
       \end{lem}
}

       \proof
       Let $\varphi$ be a smooth function such that $\min_{[0,T)\times\Rpp\times(0,1)}((V+\xi f)-\varphi)(t,x,p)=((V+\xi f)-\varphi)(t_0,x_0,p_0)=0$, $\xi>0$.
         Since $ f$ is a smooth function, the function $\psi\defeq\varphi-\xi f$ is a test function for $V$ at $(t_0,x_0,p_0)$.
         \\
        We consider $\eta \in \SmD$, recall Remark \ref{re opcont}. Using the definition of $\mathrm H\mage{^\eta}$, the inequality
        \begin{align*}
       & \hat\mu_{Y}(\cdot,(\psi+\xi f)(\cdot,p_0),\mathrm{D}(\psi+\xi f)(\cdot,p_0),\pr{\eta})(t_0,x_0)\\&\quad\ge \hat\mu_{Y}(\cdot,\psi(\cdot,p_0),\mathrm{D}\psi(\cdot,p_0),\pr{\eta})(t_0,x_0)\\&
        \quad\quad+\hat\mu_{Y}(\cdot,\psi(\cdot,p_0),\mathrm{D}(\psi+\xi f)(\cdot,p_0),\pr{\eta})(t_0,x_0)-\hat\mu_{Y}(\cdot,\psi(\cdot,p_0),\mathrm{D}\psi(\cdot,p_0),\pr{\eta})(t_0,x_0)\,,
        \end{align*}
        (recall Remark \ref{re ope increasing}) and the assumption on $\hat \mu_Y$, we obtain
         \begin{align*}
         {\mathrm{H}}^{\eta}\varphi(t_0,x_0,p_0)\ge {\mathrm{H}}^{\eta}\psi(t_0,x_0,p_0) + \mathfrak{A} + \mathfrak{B}\,,
         \end{align*}
         where
         \begin{align*}
         \mathfrak{A} &= \mage{\xi} (\eta^1)^2\left( -\partial_t \phi - \mage{\bar c} |  \mathrm D_p \phi| |\pr{\eta}| \mage{- \frac12 \|\pr{\eta}\|^2 \mathrm D_{pp} \phi} \right)\mage{(t_0,x_0,p_0)}\,,
         \\
         \mathfrak{B} &= \mage{\xi}(\eta^1)^2 \left( -\partial_t h - L | \diag{\cdot} \mathrm{D}_x h| - |\mage{\mu^\top_X}(\cdot)\mathrm{D}_x h| - \frac12|\Tr\left[\sigma_X\sigma_X^\top(\cdot)\mathrm{D}_{xx}\mage{h}\right] |\right)\mage{(t_0,x_0,p_0)}\,.
          \end{align*}
          We will now  give a lower bound for both terms.\\
         1. For the first term we compute, observing that $\partial_t \phi = - \kappa \phi$, $\mathrm D_p\phi(t,p) = 2 \bar{c} e^{\kappa(T-t)-4\bar{c}p}$  and $\mathrm D_{pp}\phi(t,p)=-8\bar{c}^2e^{\kappa(T-t)-4\bar{c}p}$,
          \begin{align*}
          \mathfrak{A} &=  \xi (\eta^1)^2  e^{\kappa(T-t_{\mage{0}})} \left( \kappa(\theta - \mage{\frac{e^{-4\bar{c}p_{\mage{0}}}}{2}})  - \mage 2\bar{c}^2 e^{-\mage 4\bar{c}p_{\mage{0}}} |\pr{\eta}| \mage{+ 4\bar{c}^2 e^{-4\bar{c}p_{\mage{0}}} \|\pr{\eta}\|^2} \right)
          \\
          & \ge  \xi (\eta^1)^2  e^{\kappa(T-t_{\mage{0}})-\mage 4\bar{c}p_{\mage{0}}} \left( (\theta - 1) \mage{- 2\bar{c}^2(1+\|\pr{\eta}\|^2) + 4 \bar{c}^2 \|\pr{\eta}\|^2} \right)\,,
          \end{align*}
          where we used the fact that $\kappa, \theta \mage > 1$. 
          \\
          Setting \mage{$\theta\defeq 4\bar c^2+1$},
          we obtain observing that $(\eta^1)^2(1+\|\pr{\eta}\|^2)=\|\eta\|^2=1$,
          \begin{align}\label{eq majo term A}
           \mathfrak{A} \ge \mage{2{\xi}\bar c^2e^{-4\bar{c}}}\,.
          \end{align}
  2. For the second term, observing that
  \begin{gather*}
  \partial_t h = -\kappa h, \quad \mathrm{D}_x h(t,x) = e^{\kappa(T-t_{\mage{0}})} (2k|x|^{2k-1}-2|x|^{-3})\1  \\ \text{and}\\
   \mathrm{D}_{xx} h(t,x) = e^{\kappa(T-t_{\mage{0}})} (2k(2k-1)|x|^{2k-2}+6|x|^{-4})\textbf{I} \,,
   \end{gather*}
   we compute,
 \begin{align*}
   \mathfrak{B} &\ge \xi (\eta^1)^2 e^{\kappa(T-t_0)} \left( \kappa
   -2(k+1)(L+\Lambda) - \Lambda^2\frac d2(2k(2k-1)+6)\right) \left(|x|^{2k}+ |x|^{-2} \right)
   \ge \xi(\eta^1)^2\,,
  \end{align*}
with $d$ the dimension of $X$ and for $\kappa$ large enough. In particular, we get
$
 \mathfrak{B} \ge 0\,.
$
\\
Combining this last inequality with \eqref{eq majo term A}, we obtain
  \begin{align*}
         {\mathrm{H}}^{\eta}\varphi(t_0,x_0,p_0)\ge {\mathrm{H}}^{\eta}\psi(t_0,x_0,p_0) + \mage{2{\xi}\bar c^2e^{-4\bar{c}}}\,.
         \end{align*}
3. We thus get that
\begin{align*}
 \sup_{\eta \in \SmD}{\mathrm{H}}^{\eta}\varphi(t_0,x_0,p_0)\ge \sup_{\eta \in \SmD}{\mathrm{H}}^{\eta}\psi(t_0,x_0,p_0) +\mage{2{\xi}\bar c^2e^{-4\bar{c}}}\,.
\end{align*}
We also observe that $V+\xi f \ge V+\mage{4\xi\bar{c}^2}$ with the above choice of $\theta$.
The proof is concluded by using the super-solution property of $\psi$, recall Remark \ref{re opcont}.
\eproof

 \bleu
 {
     \begin{lem}[Modulus of continuity]\label{le ModulusCont0} Let $(b,x,r,p,q)\in \R\times (\Rpp)^2\times\bleu{[0,1]}^2$ and $(y',y) \in \R^2$ with $y'>y$. Moreover, for $\varepsilon>0$,   let $\cX$ and $\cR\in\s^{d+1}\times\s^{d+1}$ being such that
         \begin{equation}\label{A10}
           \left(\begin{matrix} \cX&0\\0&-\cR \end{matrix}\right)\leq
           \frac3{\varepsilon}\left(\begin{matrix} I &-I\\-I& I \end{matrix}\right)\,.
         \end{equation}
           Setting
         $\delta= \frac{2}\varepsilon \left(\begin{matrix}x-r\\p-q\end{matrix} \right)$,
         $\Theta=(t,r,y,b,\delta,\cR)$, $\Theta'=(t,x,y',b,\delta,\cX)$ then
          \begin{align*}
           &\sup_{\eta\in\cS}{\mathrm  H}^\eta(\Theta)- \sup_{\eta\in\cS}{\mathrm  H}^\eta(\Theta')
            \leq C\mage{\left(|x-r|+\frac{1}\varepsilon|x-r|^2\mage{(1+|x|+|r|)}
		+\frac{1}\varepsilon|p-q|^2\right)}\,,
         \end{align*}
         for some constant $C>0$.
       \end{lem}
}
       \proof
       Consider $\Theta$ and $\Theta'$ defined in the theorem. We observe that
         \begin{align*}
          &\sup_{\eta\in\cS}{\mathrm { H}}^\eta(\Theta)- \sup_{\eta\in\cS}{\mathrm { H}}^\eta(\Theta')
           \le \sup_{\eta\in\SmD}\left\{{\mathrm { H}}^\eta(\Theta)
           -{\mathrm { H}}^\eta(\Theta')\right\}\;,
           \end{align*}
           (recall Remark \ref{re opcont}).\\
         For $\eta \in \SmD$, using the definition of $\mathrm { H}\mage{^\eta}$ and the Lipschitz property of $\mu_X$, we then compute
         \begin{align*}
          {\mathrm { H}}^\eta(\Theta)
           -{\mathrm { H}}^\eta(\Theta') \le (\eta^1)^2\left( \mathfrak{A}
           +\frac{C}{\varepsilon}|x-r|^2 + \mathfrak{B} \right)\,,
         \end{align*}
         with \mage{$C>0$} and
         \begin{align*}
          \mathfrak{A} &= \hat{\mu}_Y(t,r,y,\delta,\pr{\eta})-\hat{\mu}_Y(t,x,y',\delta,\pr{\eta})
          \text{ and }
          \mathfrak{B} = -\frac12 \Tr\left[\bar{\sigma}\bar{\sigma}^\top(r,\pr{\eta})\cR \right ] + \frac12 \Tr \left [\bar{\sigma}\bar{\sigma}^\top(x,\pr{\eta})\cX \right ]\,.
         \end{align*}
	  Since $y'>y$, using the monotonicity property of $\hat{\mu}_Y$ (recall Remark \ref{re ope increasing}) we have
	  \begin{align*}
	   \mathfrak{A} &\le |\hat{\mu}_Y(t,r,y',\delta,\pr{\eta})-\hat{\mu}_Y(t,x,y',\delta,\pr{\eta})|
	   \\
	   &\le C\left(|x-r|+\frac{1}\varepsilon|x-r|^2\mage{(1+|x|+|r|)+\Lambda\frac{|\pr{\eta}|}{\varepsilon}|x-r||p-q|}
		+\frac{|\pr{\eta}|}\varepsilon|p-q||\sigma^{-1}(r)-\sigma^{-1}(x)|
		\right)\,,
	  \end{align*}
	  where for the last inequality we used the definition of $\hat{\mu}_Y$, the Lipschitz property of $\mu_Y$ \mage{and the bound $\Lambda$ of $\sigma^{-1}$ (recall Standing Assumption)}. Using then Young's inequality and the Lipschitz continuity of $\sigma^{-1}$, we finally obtain
	  \begin{align*}
	   \mathfrak{A} &\le C(1 \mage{+ |\pr{\eta}|})\left(|x-r|+\frac{1}\varepsilon|x-r|^2\mage{(1+|x|+|r|)}
		+\frac{1}\varepsilon|p-q|^2\right)\;.
	  \end{align*}
	  For the second order term $\mathfrak{B}$, we only have to use \cite[Example 3.6]{CIL92}, especially equation \mage{\eqref{A10}}, recalling that $x \mapsto \bar{\sigma}(x,a)$ has the same Lipschitz constant as $\sigma_X$ by construction. We thus obtain
	  \begin{align*}
	   \mathfrak{B} \le \frac{C}\varepsilon|x-r|^2\;.
	  \end{align*}
        We thus have
        \begin{align*}
         {\mathrm { H}}^\eta(\Theta)
           -{\mathrm { H}}^\eta(\Theta')
           \le C(\eta^1)^2(1+\mage{ |\pr{\eta}|})\mage{\left(|x-r|+\frac{1}\varepsilon|x-r|^2\mage{(1+|x|+|r|)}
		+\frac{1}\varepsilon|p-q|^2\right)}\;.
        \end{align*}
        The proof is concluded by observing that $|\eta^1|\le 1$ and $|\eta^1\pr{\eta}| \le 1$.
        \eproof

        \subsection{The comparison principle}\label{se CPR}

        We can now prove the comparison principle for \eqref{eq TargetSolAlt}.

       \begin{theo}\label{CparisonPrinciple}
        Let ${V}$ (resp. ${U}$) be a non-negative lower semi-continuous (resp. upper semi-continuous) map satisfying a polynomial growth of order $k\ge1$ (defined in Assumption \ref{assump polygrowth0}) on $[0,T]\times\Rpp\times\bleu{[0,1]}$. Moreover assume that,
        \begin{itemize}[noitemsep,nolistsep]
        \item  on $[0,T)\times\Rpp\times(0,1)$, ${U}$ is a viscosity sub-solution of \eqref{eq TargetSolAlt} and ${V}$ is a viscosity super-solution of \eqref{eq TargetSolAlt},
        \item for all $(t,x)\in[0,T)\times\Rpp$, ${V}(\cdot,\mage{0})\ge{U}(\cdot,0)$ and ${V}(\cdot,1)\ge{U}(\cdot,1)$,
        \item for all $(x,p)\in\Rpp\times\bleu{[0,1]}$, ${V}(T,\cdot)\ge {U}(T,\cdot)$,
        \end{itemize}
         then ${V}\geq{U}$ on $[0,T]\times\Rpp\times[0,1]$.
       \end{theo}

\begin{rem}\label{re bord}
The boundary conditions are necessary for the comparison to hold. Indeed let us assume that $\mu_Y\equiv 0$ and consider on $[0,T]\times[0,1]$ the function
$$v_{\lambda,\gamma}(t,p)\defeq 2T\kappa+(t-T)\left[\lambda p+\gamma(1-p)\right],\,(\lambda,\gamma)\in[0,\kappa]\times[0,\kappa]\,,$$
for some $\kappa>0$. Note that $v_{\lambda,\gamma}(T,\cdot)=2T\kappa$ and $v_{\lambda,\gamma}\ge 0$ for all $(\lambda,\gamma)\ge[0,\kappa]\times[0,\kappa]$. We can easily prove that all the functions that belong to the family $(v_{\lambda,\gamma}(t,p))_{\lambda,\gamma}$ are a viscosity sub-solution of \eqref{eq subsolBET} and a viscosity super-solution of \eqref{eq supersolBET}. There is thus no unique solution to the system \eqref{eq subsolBET}-\eqref{eq supersolBET}.
However, if we set the following boundary conditions
$$v_{\lambda,\gamma}(\cdot,0)=2T\kappa\quad\text{and}\quad v_{\lambda,\gamma}(\cdot,1)=2T\kappa\,,$$
we thus obtain that $\lambda=\gamma=0$ and that $v_{0,0}=2T\kappa$ is the unique solution to the above system.
\end{rem}

{\begin{cor}
         Assume that $v^*(\cdot,0)\le v_*(\cdot,0)$ and $v^*(\cdot,1)\le v_*(\cdot,1)$ on $[0,T)\times\Rpp$ and assume that $v^*(T,\cdot)\le v_*(T,\cdot)$ on $\Rpp\times[0,1]$. Then the function $v$ is continuous on $[0,T]\times\Rpp\times[0,1]$ and is the unique  viscosity solution, in the class of function
         with polynomial growth, of
         \begin{align*}
         \mathcal{H}v =0\; \text{ on } [0,T)\times\Rpp\times(0,1)\,.
         \end{align*}
\end{cor}}

\noindent{\bf Proof of Theorem \ref{CparisonPrinciple}}
          Let us now define on $(0,T]\times\Rpp\times[0,1]$ the following non-negative auxiliary function
         \begin{align*}
         {{V}}_{\xi}(t,x,p)&\defeq ({V}+\xi f)(t,x,p) + \frac{\xi}{t}\,,
         \end{align*}
         with $f$ defined in Lemma \ref{le StrictSupersol}. Using Lemma \ref{le StrictSupersol}, it is easily seen that $ {{V}}_{\xi}$ is a strict super-solution
         of \eqref{eq TargetSolAlt}, namely it satisfies \eqref{eq strict supersol}.
         \\
         We also introduce $U_\xi(t,x,p)= U(t,x,p) - \xi h(t,x)$, recall Lemma \ref{le StrictSupersol}. By an easy adaptation of the proof of Lemma \ref{le StrictSupersol}, we have that $U_\xi$ is still a sub-solution to \eqref{eq subsolBET}.

         \mage{We will show that ${U}-{V}\leq0$ on $\mage(0,T]\times\Rpp\times\bleu{[0,1]}$.} To this aim we will first prove by contradiction that for all $\xi>0$ we have ${{U}}_{\xi}-{{V}}_{\xi }\leq0$ and the result will follow sending $\xi$ to zero.\\
         \textbf{Step 1.}
         We assume to the contrary that \mage{there exists $\xi>0$ such that $$\sup_{(0,T]\times\Rpp\times[0,1]}({U}_{\xi}-{V}_{\xi})(t,x,p)=({U}_{\xi}-{V}_{\xi})(\hat{t},\hat{x},\hat{p}) =\gamma>0\,.$$}
      \mage{Observe that as $V_{\xi }>0$ the previous expression implies that $$U_{\xi }(\hat t, \hat x,\hat p)> 0\,.$$}
         For $\varepsilon>0$, we define on $(0,T]\times(\Rpp)^2\times\bleu{[0,1]}^2$
        \begin{align*}
          \Psi_{\varepsilon}(t,x,r,p,q) \defeq {U}_\xi(t,x,p)-{V}_{\xi}(t,r,q)-\frac{1}{\varepsilon}\left(|x-r|^2 +|p-q|^2\right)\,. 
        \end{align*}
       Using the growth conditions and semi-continuity of  ${U}$ and ${V}$, it follows that for $\varepsilon>0$ the function $\Psi_{\varepsilon}$ admits a maximum $M_\varepsilon$ at $(t_{\varepsilon},x_{\varepsilon},r_{\varepsilon},p_{\varepsilon},
       q_{\varepsilon})$ on $\mage(0,T]\times{\Rpp}^2\times\bleu{[0,1]}^2$. Moreover the inequality $\Psi_{\varepsilon}(t_{\varepsilon},x_{\varepsilon},r_{\varepsilon},p_{\varepsilon},
       q_{\varepsilon})\ge\Psi_{\varepsilon}(\hat{t},\hat{x},\hat{x},\hat{p},\hat{p}) = \gamma$ combined with the growth condition on $U$, $V$ and the definition of $f$ and $h$, implies that $t_\varepsilon$, $x_{\varepsilon}$ and $r_{\varepsilon}$ are in compact set $\cT\times \cK \subset (0,T]\times \Rpp$. Let $(\hat{t},\hat{x},\hat{p})\in\cT\times \cK \times [0,1]$ be a limit point of $(t_{\varepsilon}, x_{\varepsilon},p_{\varepsilon})$. Using \cite[Lemma 3.1]{CIL92}
       we obtain that
       \begin{align}
       \begin{cases}\label{eq cv}
         \lim_{\varepsilon \downarrow 0} \frac{1}{\varepsilon}
	      \left(|x_\varepsilon-r_\varepsilon|^2 +|p_\varepsilon-q_\varepsilon|^2\right) &= 0\,,
	      \\
	 \lim_{\varepsilon \downarrow 0} M_\varepsilon = ({U}_{\xi}-{V}_{\xi})(\hat{t},\hat{x},\hat{p}) \,.
       \end{cases}
       \end{align}
                If $(\hat{t},\hat{x},\hat{p}) \in \mage (0,T)\times\Rpp\times\set{0,1}$ or $(\hat{t},\hat{x},\hat{p})\in \{T\}\times\Rpp\times\bleu{[0,1]}$, the
         assumptions on ${V}$ and ${U}$ on these boundaries of the domain lead to a contradiction.
        \\
          We thus now assume that $0<\hat{t}<T$ and $0<\hat{p}<1$. In particular, up to a subsequence,
          \begin{align}\label{eq intermCP0}
       U_{\xi }(t_\varepsilon,x_\varepsilon,p_\varepsilon)> 0\,.
       \end{align}
       \\
       \textbf{Step 2.}
       From Ishii's Lemma (see \cite[Theorem 8.3]{CIL92}) we get the existence of real coefficients $b^1_{\varepsilon},b^2_{\varepsilon}$, a vector $d_{\varepsilon}$ and two symmetric matrices $\mathcal{X}_{\varepsilon}$ and $\mathcal{R}_{\varepsilon}$ being such that
       \begin{align*}
         (b^1_{\varepsilon},d_{\varepsilon},\mathcal{X}_{\varepsilon})
         \in\bar{\mathcal{J}}^+_{\bar{\mathcal{O}}}{U}_{\mage\xi}(t_{\varepsilon},x_{\varepsilon},p_{\varepsilon}) &\quad\text{and}\quad(-b^2_{\varepsilon},d_{\varepsilon},\mathcal{R}_{\varepsilon})
         \in\bar{\mathcal{J}}^-_{\bar{\mathcal{O}}}{V}_{\xi}(t_{\varepsilon},r_{\varepsilon},q_{\varepsilon})\,,
       \end{align*}
       with {$\bar{\mathcal{O}}\defeq\mage{(0,T)\times \Rpp }\times(0,1)$} and $\bar{\mathcal{J}}^+$ (resp. $\bar{\mathcal{J}}^-$) the limiting second order super-jet (resp. sub-jet) of ${U}_{\mage\xi}$ (resp. ${V}_{\xi}$) at $(t_{\varepsilon},x_{\varepsilon},p_{\varepsilon}) \in\bar{\mathcal{O}}$ (resp. $(t_{\varepsilon},r_{\varepsilon},q_{\varepsilon}) \in\bar{\mathcal{O}}$) and
       where
       \begin{subequations}\label{eq: Der}
       \begin{align}
         b^1_{\varepsilon}+b^2_{\varepsilon}&\defeq 0 \,,\label{eq: Der b1}\\
        d_{\varepsilon}
         &\defeq\frac{2}{\varepsilon}\left(\begin{matrix}
                             x_{\varepsilon}-r_{\varepsilon}
                              \\
                              p_{\varepsilon}-q_{\varepsilon}
                              \end{matrix}\right)\,,\label{eq: Der e}
        \\
        \left(\begin{matrix}
               \mathcal{X}_{\varepsilon} & 0 \\
               0 & - \mathcal{R}_{\varepsilon}
              \end{matrix}
\right) &\le \frac3{\epsilon}
\left(\begin{matrix}
               I  & -I \\
               -I &  I
              \end{matrix}
\right)\,.
        \end{align}
       \end{subequations}
       \mage{It follows from the definition of ${U}_\xi$ and ${V}_{\xi}$ that they are respectively sub-/super-solution of \eqref{eq TargetSolAlt}. As a consequence using \eqref{eq intermCP0} and Lemma \ref{le StrictSupersol} we obtain}
       \begin{align*}
         \sup_{\eta\in\cS}{\mathrm H}^{\eta}\left(t_{\varepsilon},x_{\varepsilon},{U}_{\mage\xi}(t_{\varepsilon},x_{\varepsilon},p_{\varepsilon}),
         b^1_{\varepsilon},d_{\varepsilon},\mathcal{X}_{\varepsilon}\right)&\le 0\,, \\
         \sup_{\eta\in\cS}{\mathrm H}^{\eta}\left(t_{\varepsilon},r_{\varepsilon},{V}_{\xi}(t_{\varepsilon},r_{\varepsilon},q_{\varepsilon}),
         -b^2_{\varepsilon},d_{\varepsilon},\mathcal{R}_{\varepsilon}\right)&\ge\xi \rho >0\,.
       \end{align*}
        Hence
       \begin{gather}
          \sup_{\eta\in\cS}{\mathrm H}^{\eta}\left(t_{\varepsilon},r_{\varepsilon},{V}_{\xi}(t_{\varepsilon},r_{\varepsilon},q_{\varepsilon}),
         -b^2_{\varepsilon},d_{\varepsilon},\mathcal{R}_{\varepsilon}\right) - \sup_{\eta\in\cS}{\mathrm H}^{\eta}\left(t_{\varepsilon},x_{\varepsilon},{U}_{\mage\xi}(t_{\varepsilon},x_{\varepsilon},p_{\varepsilon}),
         b^1_{\varepsilon},d_{\varepsilon},\mathcal{X}_{\varepsilon}\right)
         \ge \xi \rho \,.\label{eq CtradictoryModCont}
       \end{gather}
       \textbf{Step 3.}
       On the other hand we know from Lemma \ref{le ModulusCont0} and \eqref{eq: Der} that \mage{there exists $C>0$ such that}
       \begin{gather*}
         \sup_{\eta\in\cS}{\mathrm H}^{\eta}\left(t_{\varepsilon},r_{\varepsilon},{V}_{\xi}(t_{\varepsilon},r_{\varepsilon},q_{\varepsilon}),
         -b^2_{\varepsilon},d_{\varepsilon},\mathcal{R}_{\varepsilon}\right)- \sup_{\eta\in\cS}{\mathrm H}^{\eta}\left(t_{\varepsilon},x_{\varepsilon},{U}_{\mage\xi}(t_{\varepsilon},x_{\varepsilon},p_{\varepsilon}),
         b^1_{\varepsilon},d_{\varepsilon},\mathcal{X}_{\varepsilon}\right)
         \nonumber\\
         \le C\left(|x_\varepsilon-r_\varepsilon|+\frac{1}\varepsilon|x_\varepsilon-r_\varepsilon|^2\mage{(1+|x_\varepsilon|+|r_\varepsilon|)}
		+\frac{1}\varepsilon|p_\varepsilon-q_\varepsilon|^2\right)\label{eq imtermCP1}
         \,.
       \end{gather*}
       Now sending  $\varepsilon$  to zero and using \eqref{eq cv} we obtain
       that the last inequality is non positive. Thus we obtain a contradiction to \eqref{eq CtradictoryModCont} and ${U}_{\mage\xi}\leq{V}_{\xi}$ for all $\xi>0$ on $\mage(0,T]\times\Rpp\times\bleu{[0,1]}$. This gives the required result by sending $\xi$ to zero.
      \qed

  \section{Application to the quantile hedging of Bermudan options}\label{se appli}
In this section we are interested in the case where the loss function is the indicator function leading to a quantile hedging problem. More precisely we are interested in the Bermudan version of this problem, i.e. we define for $(t,x,p)\in[0,T]\times\Rpp\times[0,1]$
\begin{eqnarray*}
&v(t,x,p):=\inf\left\{y\ge 0\,:\,\exists\, \nu \in \mathcal{U}_{t,x,y}\mbox{ s.t. }\p\left[\bigcap_{s\in \bT_{t}} \mathrm{S}^{t,x,y,\nu}_{s}\right]\ge p\right\}\;,&\nonumber \\
&\mbox{ with } \;
 \mathrm S^{t,x,y,\nu}_{s}:=\left\{
 \begin{array}{lcl}
 \Omega & \mbox{if }&    s \le t \\
 \{Y^{t,x,y,\nu}_{s}\ge g(s,X^{t,x}_{s})\}   & \mbox{if }& s>t
\end{array}\right.\,,
\nonumber\\&\mbox{ and } \; \mathbb{T}_{t}\defeq \set{t_0=0\leq \dots \leq t_i\leq \dots\leq t_n=T}\cap (t,T]
\label{eq applidef}\,,
\end{eqnarray*}
    where $g: [0,T]\times\Rpp\rightarrow\R^+$ and $x\in\Rpp\mapsto g(t,x)$ is Lipschitz continuous for all $t\in[0,T]$ and where $\mathcal{U}_{t,x,y}$ is defined in Section \ref{se pbstatement}.
 Observe that $v(t,\cdot)$ must be  interpreted as a
{continuation}  value, i.e.~the price at time $t$ knowing that the option
has not been exercised on $[0,t]$. In particular,   $v(T,\cdot)=0$ on $\Rpp\times[0,1]$.

 This problem is equivalent to (see \eqref{eq Pb1Alt})
       \begin{equation*}
         v(t,x,p) = \inf \left\{\begin{matrix}y\in\R^+: \exists \hspace{1mm} (\nu,\alpha)\in \mathcal{U}_{t,x,y}\times\mathcal{A}_{t,p}\hspace{2mm}\text{s.t.}\,Y^{t,x,y,\nu}_s\geq g(s,X^{t,x}_s)\1_{\{P^{t,p,\alpha}_s>0\}},\,\forall\,s\in\mathbb{T}_{t}\end{matrix}\right\}\,,
       \end{equation*}
       where for $p\in[0,1]$, $\mathcal{A}_{t,p}$ is the set of $\R^d$-valued $\F$-progressively measurable and square integrable processes $\alpha$ such that $P^{t,p,\alpha}_\cdot\in[0,1]$ on $[t,T]$. 

\bleu{
The aim of this section is to give a  characterisation of $v$ as the unique solution to a sequence of PDEs. To this end, and in view of the previous section \mage{(recall Remark \ref{re bord})}, we need in particular the knowledge of $v$ on the boundary $p=0$ and $p=1$. We will show that,
as expected, $v(\cdot,0)=0$ and $v(\cdot,1)=\bar{v}(\cdot)$, where $\bar{v}$ is the super-replication price of the Bermudan {o}ption with exercise
price $g$. Moreover, we  assume that Assumption \ref{assump polygrowth0} holds for $\bar{v}$ in this Bermudan setting. Therefore, Remark \ref{re polygrowth0} is still
valid for $v$. Precisely \mage{from standard results in stochastic control theory (see \cite{P09} for instance) and from \cite{BET09}},
we have the following characterisation of $\bar{v}$.
 \begin{prop}\label{pr carac vbar}
         Fix $1\le i\le n$.
         \\ The function $\bar{v}$ is continuous on $[t_{i-1},t_i)\times\Rpp$ and is the unique viscosity solution, in the class of function with polynomial growth, of
            \begin{align*}
               \min\left\{\bar v(t,x),\left\{\begin{matrix}-\partial_t\varphi(t,x)+\mu_Y(t,x,y,\diag{x}\mathrm D_{x}\varphi(t,x))-\mu_X^\top(x)\mathrm D_{x}\varphi(t,x)\\-\frac{1}{2}\Tr\left[\sigma_X\sigma_X^\top(x)\mathrm D_{xx}\varphi(t,x)\right]\end{matrix}\right\}\right\}=0\,,
            \end{align*}
            with terminal condition at time $t_i$,
            $$
            \lim_{t \uparrow t_i} v(t_{i},\cdot) = \bar v(t_{i},\cdot)\vee g(t_i,\cdot)\;.
            $$
      \end{prop}
}
We can now state the main result of the section which is the full PDE characterisation of $v$, the quantile hedging price of the Bermudan {o}ption with exercise price $g$.

\begin{theo} \label{th objectif} Fix $1\le i\le n$.
         \\ The function $v$ is continuous on $[t_{i-1},t_i)\times\Rpp\times[0,1]$ and is the unique  viscosity solution, in the class of functions
         with polynomial growth, of
         \begin{align*}
         \mathcal{H}v(t,x,p) =0\;, \text{ for } (t,x,p) \in [t_{i-1},t_i)\times\Rpp\times(0,1)\,,
         \end{align*}
         with boundary conditions $v(t,x,0)=0$, $v(t,x,1)=\bar{v}(t,x)$, $(t,x)\in\mage{[t_{i-1},t_{i})}\times\Rpp$ and
         \begin{align*}
         \lim_{t \uparrow t_i} v(t,x,p) = \overline{\text{conv}}\left(v(t_{i},x,p)\vee g(t_{i},x)\1_{\set{p>0}}\right)\;,\quad (x,p) \in \Rpp\times [0,1].
         \end{align*}
\end{theo}

Using \cite[Proposition 3.3 (a)]{BBC14}, we observe that the terminal condition at  time $t_i$, $1 \le i \le n$, can be easily computed. More precisely it is obtained by
applying \cite[Lemma 3.1(a)]{BBC14}, the fact that $v(t_{i+1},\cdot,0)=0$ and the definition of $p_g$.

\begin{rem} \label{re facelift} (i) The boundary condition at time $t_i$, $1 \le i \le n$, and for $(x,p)\in\Rpp\times[0,1]$, is given by
\begin{align}\label{eq facelift}
\overline{\text{conv}}\left(v(t_{i},x,p)\vee g(t_{i},x)\1_{\set{p>0}}\right) = v(t_i,x,p)\vee \tilde{g}(t_i,x,p)\,,
\end{align}
where for $(t,x,p) \in
[0,T]\times\Rpp\times\R$, $\tilde{g}$ is the following `facelift' of $g$
 \begin{align*}
     \tilde{g}(t,x,p) = q_g(t,x) p\1_{\{0 \le p \le 1 \}}
+\infty\1_{\{p > 1\}}  \,,
    \end{align*}
    with
    \begin{align*}
    q_g(t,x) := \frac{g(t,x)}{p_g(t,x)} \1_{\set{p_g(t,x)>0}} \;
\mbox{ and }\;
      p_g(t,x) := \sup \set{p \in \R \,|\, v(t,x,p)=g(t,x) }\;
\wedge 1\,.
    \end{align*}
    (ii) In particular, at time $T$, the terminal condition is given by
    $(t,x,p) \mapsto p g(t,x)$, which was already observed in \cite[Proposition 3.2]{BET09}.
\end{rem}

\subsection{Proof of Theorem \ref{th objectif}}

We now turn to the proof of the main result of this section. As usual in the case of Bermudan option,
the proof is done by induction on the time interval $[t_i,t_{i+1})$, $0 \le i \le n-1$. The main difficulty here is
the characterisation of $v$ on the boundaries of the domain, specially the time-boundary for which a facelifting phenomenon appears.


The results stated in this section are a direct consequence of the geometric dynamic programming principle, see \cite{BET09, ST02GEO,ST02}. In our framework, we obtain from \cite[Theorem 2.1]{BV10}, the following geometric dynamic programming principle,\\
       \newline
         \textbf{(GDP1)} Fix $1\le i\le n$ and $(t,x,p)\in[t_{i-1},t_i)\times\Rpp\times[0,1]$. If $y>v(t,x,p)$, then there exists a $(\nu,\alpha)\in\mathcal{U}_{t,x,y}\times\mathcal{A}_{t,p}$ such that for all stopping times $\theta\le t_i$
       \begin{align*}
         Y^{t,x,y,\nu}_\theta\geq v\left(\theta,X^{t,x}_\theta,P^{t,p,\alpha}_\theta\right)\1_{\set{\theta < t_i}}
         + (v\vee g)\left(t_i,X^{t,x}_{t_i},P^{t,p,\alpha}_{t_i}\right)\1_{\set{\theta = t_i}}\,.
        \end{align*}
         \textbf{(GDP2)} Fix $1\le i\le n$ and $(t,x,p)\in[t_{i-1},t_i)\times\Rpp\times[0,1]$. For every $y<v(t,x,p)$, $(\nu,\alpha)\in\mathcal{U}_{t,x,y}\times\mathcal{A}_{t,p}$ and all stopping times $\theta\le t_i$
       \begin{align*}
         \p\left[Y^{t,x,y,\nu}_\theta > v\left(\theta,X^{t,x}_\theta,P^{t,p,\alpha}_\theta\right)\1_{\set{\theta < t_i}}
         + (v\vee g)\left(t_i,X^{t,x}_{t_i},P^{t,p,\alpha}_{t_i}\right)\1_{\set{\theta = t_i}}\right] < 1\,,
        \end{align*}
with the notation $g(t_i,x,p)\defeq g(t_i,x)\1_{p>0}+\infty\1_{p>1},\,1\le i\le n$ and $x\in\Rpp$.

\vspace{5pt}
\noindent
Let us now proceed with the proof of  Theorem \ref{th objectif}.\\
For $i\le n-1$, we now assume that $v(t_{i+1},\cdot)$ is continuous, $v(t_{i+1},\cdot,0)=0$ and $v(t_{i+1},\cdot,1)=\bar v(t_{i+1},\cdot)$. (Observe that this is the case by convention at time $T$ as $v(T,\cdot)=0=\bar v(T,\cdot)$).
\\
To clarify the arguments, we introduce the following function on \mage{$[t_i,t_{i+1}]\times\Rpp\times[0,1]$}
       \begin{align*}\hat v(t,x,p)\defeq \begin{cases} v(t,x,p)\,\text{on}\, [t_{i},t_{i+1})\times\Rpp\times[0,1]\\
                                        v(t_{i+1},x,p)\vee g(t_{i+1},x)\1_{p>0}\,\text{on}\, \Rpp\times[0,1]
       \end{cases}\,.
       \end{align*}
\textbf{Step 1} Characterisation on $[t_i,t_{i+1})\times\Rpp\times[0,1]$.\\
From $\mathbf{(GDP1)-(GDP2)}$, combining the results of \cite[Theorem 2.1]{BET09} and Theorem \ref{th AltPDECharac0}, we obtain that
$\hat{v}$ is a viscosity solution of
  \begin{align*}
         \mathcal{H}\hat{v} =0\;, \text{ on }  [t_i,t_{i+1})\times\Rpp\times(0,1)\,.
         \end{align*}
Moreover applying \cite[Theorem 3.1]{BET09} we obtain that, on $[t_i,t_{i+1})\times\Rpp$,
\begin{align} \label{eq p-boundaries}
\hat{v}^*(\cdot,1)=\hat{v}_*(\cdot,1)=\bar{v}(\cdot) \text{  and  } \hat{v}^*(\cdot,0)=\hat{v}_*(\cdot,0)=0\,.
\end{align}
\\
\textbf{Step 2} Characterisation on $\set{t_{i+1}}\times\Rpp\times[0,1]$.\\
\textbf{Step 2.a} We first prove that
\begin{align*}
\hat{v}^*(t_{i+1},x,p) \leq \overline{\text{conv}}\left(v(t_{i+1},x,p)\vee g(t_{i+1},x)\1_{\set{p>0}}\right)\;.
\end{align*}
Proceeding as in \cite[Section 5.4]{BET09} we first obtain
 \begin{align} \label{eq majo scs-v}
 \hat{v}^*(t_{i+1},x,p)\le\left(v(t_{i+1},x,p)\vee g(t_{i+1},x)\1_{\{p>0\}}\right)^*\,.
 \end{align}
Now, it follows from the sub-solution property  that, for any test function $\varphi$ such that $\max_{[\mage{t_{i},t_{i+1}})\times\Rpp\times(0,1)}(\hat{v}^*-\varphi)(t,x,p)=0$, we have $\mathrm{D}_{pp}\varphi(t,x,p)\ge 0$ \mage{(recall Remark \ref{re conv} \rm (1))}. Applying then the same argument as in \cite[Proposition 5.2]{CPT99}
we conclude that $\hat{v}^*$ is convex inside the domain. From the upper semi-continuity of $\hat{v}^*$, we obtain the convexity
property in the $p$-variable on $[0,1]$. 
\\
Combining \eqref{eq p-boundaries} and  \eqref{eq majo scs-v}, we observe that
\begin{align}\label{eq conv}
\hat{v}^*(t_{i+1},x,p)\le v(t_{i+1},x,p)\vee g(t_{i+1},x)\1_{p>0}\,.
\end{align}
We now use \eqref{eq facelift} and \eqref{eq conv}.
Indeed, we observe that, for all $p \in [0,p_g(t_{i+1},x)]$,
\begin{align*}
\hat{v}^*(t_{i+1},x,p)\le v(t_{i+1},x,p) \vee \tilde{g}(t_{i+1},x,p)\,,
\end{align*}
since $p\mapsto\hat{v}^*(\cdot,p)$ is convex, $\hat{v}^*(t_{i+1},x,0)=0$ and $$\hat{v}^*(t_{i+1},x,p_g(t_{i+1},x)) \le v(t_{i+1},x,p_g(t_{i+1},x))=g(t_{i+1},x)\,.$$
For $p \in [p_g(t_{i+1},x),1]$,
we have that
\begin{align*}
\hat{v}^*(t_{i+1},x,p)\le v(t_{i+1},x,p)\vee g(t_{i+1},x)=v(t_{i+1},x,p)=v(t_{i+1},x,p) \vee \tilde{g}(t_{i+1},x,p)\,,
\end{align*}
which concludes the proof for this step.

\textbf{Step 2.b}
We now prove
\begin{align*}
\hat{v}_*(t_{i+1},x,p) \geq \overline{\text{conv}}\left(v(t_{i+1},x,p)\vee g(t_{i+1},x)\1_{\set{p>0}}\right)\;.
\end{align*}

%
To obtain the above result, we will use the following Lemma, whose proof is postponed at the end of this section.
\begin{lem}\label{le probaequiv}
For all sequences $(t_k,x_k,p_k,y_k,\nu^k,\alpha^k)_{k\ge 1}\in[t_{i},t_{i+1})\times\Rpp\times(0,1)\times\R^+\times\mathcal{U}_{t_k,x_k,y_k}\times\mathcal{A}_{t_k,p_k},\,\mage{0\le i\le n-1}$ such that $(t_k,x_k,p_k,y_k)\rightarrow (t_{i+1},x,p,y)\in\Rpp\times[0,1]\times\R^+$, there exists a sequence of non-negative random variables $(H^k_{t_{i+1}})_{k\ge1}$, such that
 \begin{gather}\label{eq probaequivbermu2}
\limsup_{k\rightarrow\infty}\E[H_{t_{i+1}}^kY^{t_k,x_k,y_k,\nu^k}_{t_{i+1}} ]\le y\quad\text{and}\quad\liminf_{k\rightarrow\infty}H_{t_{i+1}}^k=1\,.
 \end{gather}
\end{lem}

\vspace{5pt}
Fix $(x,p)\in\Rpp\times[0,1]$. Let $(t_k,x_k,p_k)_{k\ge 1}\in [t_{i},t_{i+1})\times\Rpp\times(0,1)$ be a sequence such that $(t_k,x_k,p_k)\rightarrow (t_{i+1},x,p)$ and $\hat{v}(t_k,x_k,p_k)\rightarrow \hat{v}_*(t_{i+1},x,p)$. Set for every $k\ge 1$, $y_k\defeq \mage{\hat v}(t_k,x_k,p_k)+k^{-1}$ so that by \textbf{(GDP1)} there exists $(\nu^k,\alpha^k)\in\mathcal{U}_{t_k,x_k,y_k}\times\mathcal{A}_{t_k,p_k}$ such that
$$Y^{t_k,x_k,y_k,\nu^k}_{t_{i+1}}\ge v(t_{i+1},X^{t_k,x_k}_{t_{i+1}},P^{t_k,p_k,\alpha^k}_{t_{i+1}})\vee g(t_{i+1},X^{t_k,x_k}_{t_{i+1}})\1_{\{P^{t_k,p_k,\alpha^k}_{t_{i+1}}>0\}}\,.$$
Now first multiply by $H^k_{t_{i+1}}$ and then take the expectation and the limit to obtain with Fatou's Lemma
\begin{align*}
&\liminf_{k\rightarrow\infty}\E\left[H^k_{t_{i+1}} Y^{t_k,x_k,y_k,\nu^k}_{t_{i+1}}\right]\\&\quad\ge \E\left[\liminf_{k\rightarrow\infty}\left(H_{t_{i+1}}^k \left(v(t_{i+1},X^{t_k,x_k}_{t_{i+1}},P^{t_k,p_k,\alpha^k}_{t_{i+1}})\vee g(t_{i+1},X^{t_k,x_k}_{t_{i+1}})\1_{\{P^{t_k,p_k,\alpha^k}_{t_{i+1}}>0\}}\right)\right)\right]\,.
\end{align*}
We then use the Lipschitz continuity of $g$, the continuity of $v$, the lower semi-continuity of $r\in[0,1]\mapsto\1_{r> 0}$, the $L^1$ convergence of $(X^{t_k,x_k}_{t_{i+1}},P^{t_k,p_k,\alpha^k}_{t_{i+1}})$ towards $(x,p)$ and \eqref{eq probaequivbermu2} to obtain
\begin{align*}
\mage{\hat v}_*(t_{i+1},x,p)\ge\liminf_{k\rightarrow\infty}\E\left[H_{t_{i+1}}^k Y^{t_k,x_k,y_k,\nu^k}_{t_{i+1}}\right]&\ge v(t_{i+1},x,p)\vee g(t_{i+1},x)\1_{\{p>0\}}
\\&\ge  \overline{\text{conv}}\left(v(t_{i+1},x,p)\vee g(t_{i+1},x)\1_{\{p>0\}}\right)\,,
\end{align*}
by definition of the closed convex envelope.
\\
\textbf{Step 3} To conclude, let us observe that by using the $p$-boundary condition on $\hat{v}^*$ and $\hat{v}_*$ in equation
\eqref{eq p-boundaries}, the time-boundary condition of the previous step and the comparison theorem proved in the
last section, we obtain that $\hat{v}$ is continuous on $[t_i,t_{i+1})\times\Rpp\times[0,1]$. \eproof

\vspace{20pt}
\textbf{Proof of Lemma \ref{le probaequiv}.}  For ease of notations, we introduce $Y^k := Y^{t_k,x_k,y_k,\nu^k}$, $X^k := X^{t_k,x_k}$.
For later use, let us also observe that under \mage{the standing assumptions on the coefficients of $X^k$}, the following holds true
\begin{align}
\label{eq borne Xk}
\E[ |X^k_t|^q] \leq C_q, \; \text{ for all } q \ge 1\;,
\end{align}
\mage{where $C_q>0$ is a constant that does not depend on $k$.}
We now define $\tilde{Y}^k := e^{-L(t-t_k)}H^k_tY^k_t$, $t \in  [t_k,t_{i+1}]$ where $H^k$ is the solution to
\begin{align*}
 H^k_t = 1 - \int_{t_k}^t L H^k_s  \left \{ \sigma^{-1}(X^k_s) \beta^k_s \right \}^\top \ud W_s\,,\quad \text{ for } t \in  [t_k,t_{i+1}]\,,
\end{align*}
with $\beta^k = (\mathrm{sign}[(\nu^k)^i])_{1 \le i \le d}$ \mage{and $L$ defined in \eqref{as coeflipY}}. \mage{As $\sigma^{-1}$ is bounded} and $|\beta^k| \le d$, we have that
\begin{align} \label{eq borne Hk}
\E[( H^k_t)^q] \leq C_q, \; \text{ for all } q \ge 1\;,
\end{align}
where $C_q>0$ is a constant that does not depend on $k$. In particular, we observe that $\liminf_{k\rightarrow\infty}H_{t_{i+1}}^k=1$.
\\
Now, applying Ito's formula, we compute
\begin{align*}
\tilde{Y}^k_t = y + \int_{t_k}^t \left( e^{-L(s-t_k)}H^k_s \mu_Y(s,X^k_s,Y^k_s,\nu^k_s) - L\tilde{Y}_s
- L e^{-L(s-t_k)}H^k_s (\nu^k_s)^\top \beta^k_s \right) \ud s + M^k_t-M^k_{t_k}\,,
\end{align*}
for some \mage{local} martingale $M^k$. It is easily seen that this local martingale is actually a super-martingale as it is bounded from below by an integrable term.
Using \eqref{as coeflipY}, observing that $(\nu^k_s)^\top \beta^k_s = |\nu^k_s| $ and recalling $Y \ge 0$, we obtain
\begin{align*}
0 \leq \E[\tilde{Y}^k_t] \leq y + C\int_{t_k}^t\left( 1 +\E[(1+|\nu^k_s|)H^k_s |X^k_s|]\right) \ud s \;.
\end{align*}
\mage{Using Cauchy-Schwarz's inequality, \eqref{eq borne Xk}-\eqref{eq borne Hk} and the square integrability of $\nu^k_s$}, we get, recalling the definition of $\tilde{Y}^k$,
\begin{align*}
\E[H^k_tY^k_t] \le e^{L(t-t_k)}\left( y + C(t-t_k) \right)\;,\; \text{ for all } t \in [t_k,t_{i+1}]\;,
\end{align*}
which concludes the proof of the lemma. \eproof
\bibliography{BibGB}

\begin{thebibliography}{10}

\bibitem{BBMZ09}
Olivier Bokanowski, Benjamin Bruder, Stefania Maroso, and Hasnaa Zidani.
\newblock Numerical approximation for a superreplication problem under gamma
  constraints.
\newblock {\em SIAM Journal on Numerical Analysis}, 47(3):2289--2320, 2009.

\bibitem{BPZ16}
Olivier Bokanowski, Athena Picarelli, and Hasnaa Zidani.
\newblock State-constrained stochastic optimal control problems via
  reachability approach.
\newblock {\em SIAM Journal on Control and Optimization}, 54(5):2568--2593,
  2016.

\bibitem{BBC14}
Bruno Bouchard, G{\'e}raldine Bouveret, and Jean-Fran{\c{c}}ois Chassagneux.
\newblock A backward dual representation for the quantile hedging of
  \uppercase{b}ermudan options.
\newblock {\em SIAM Journal on Financial Mathematics}, 7(1):215--235, 2016.

\bibitem{BD13}
Bruno Bouchard and Ngoc-Minh Dang.
\newblock Generalized stochastic target problems for pricing and partial
  hedging under loss constraints - application in optimal book liquidation.
\newblock {\em Finance and Stochastics}, 17(1):31--72, 2013.

\bibitem{bouchardbsdes}
Bruno Bouchard, Romuald Elie, and Antony R{\'e}veillac.
\newblock \uppercase{BSDE}s with weak terminal condition. preprint, 2012.
\newblock {\em Cited on}, pages 2,7.

\bibitem{BER15}
Bruno Bouchard, Romuald Elie, and Antony R{\'e}veillac.
\newblock \uppercase{BSDE}s with weak terminal condition.
\newblock {\em The Annals of Probability}, 43(2):572--604, 2015.

\bibitem{BET09}
Bruno Bouchard, Romuald Elie, and Nizar Touzi.
\newblock Stochastic target problems with controlled loss.
\newblock {\em SIAM Journal on Control and Optimization}, 48(5):3123--3150,
  2009.

\bibitem{BV12}
Bruno Bouchard and Thanh Nam~Vu.
\newblock A stochastic target approach for \uppercase{P\&L} matching problems.
\newblock {\em Mathematics of Operations Research}, 37(3):526--558, 2012.

\bibitem{BV10}
Bruno Bouchard and Thanh~Nam Vu.
\newblock The obstacle version of the geometric dynamic programming principle:
  Application to the pricing of \uppercase{a}merican options under constraints.
\newblock {\em Applied Mathematics and Optimization}, 61(2):235--265, 2010.

\bibitem{Bru05}
Benjamin Bruder.
\newblock Super-replication of \uppercase{e}uropean options with a derivative
  asset under constrained finite variation strategies.
\newblock {\em preprint}, 2005.

\bibitem{CST05}
Patrick Cheridito, H~Mete Soner, and Nizar Touzi.
\newblock The multi-dimensional super-replication problem under gamma
  constraints.
\newblock 22(5):633--666, 2005.

\bibitem{CIL92}
Michael~G Crandall, Hitoshi Ishii, and Pierre-Louis Lions.
\newblock User’s guide to viscosity solutions of second order partial
  differential equations.
\newblock {\em Bulletin of the \uppercase{a}merican Mathematical Society},
  27(1):1--67, 1992.

\bibitem{CPT99}
Jak{\v{s}}a Cvitani{\'c}, Huyen Pham, and Nizar Touzi.
\newblock A closed-form solution to the problem of super-replication under
  transaction costs.
\newblock {\em Finance and stochastics}, 3(1):35--54, 1999.

\bibitem{el1997backward}
Nicole El~Karoui, Shige Peng, and Marie~Claire Quenez.
\newblock Backward stochastic differential equations in finance.
\newblock {\em Mathematical finance}, 7(1):1--71, 1997.

\bibitem{FL99}
Hans F{\"o}llmer and Peter Leukert.
\newblock Quantile hedging.
\newblock {\em Finance and Stochastics}, 3(3):251--273, 1999.

\bibitem{IL90}
Hitoshi Ishii and Pierre-Louis Lions.
\newblock Viscosity solutions of fully nonlinear second-order elliptic partial
  differential equations.
\newblock {\em Journal of Differential equations}, 83(1):26--78, 1990.

\bibitem{JKT13}
Ying Jiao, Olivier Klopfenstein, and Peter Tankov.
\newblock Hedging under multiple risk constraints.
\newblock {\em arXiv preprint arXiv:1309.5094}, 2013.

\bibitem{M11}
Ludovic Moreau.
\newblock Stochastic target problems with controlled loss in jump diffusion
  models.
\newblock {\em SIAM Journal on Control and Optimization}, 49(6):2577--2607,
  2011.

\bibitem{P09}
Huy{\^e}n Pham.
\newblock {\em Continuous-time stochastic control and optimization with
  financial applications}, volume~61.
\newblock Springer Science \& Business Media, 2009.

\bibitem{ST02GEO}
H~Mete Soner and Nizar Touzi.
\newblock Dynamic programming for stochastic target problems and geometric
  flows.
\newblock {\em Journal of the \uppercase{e}uropean Mathematical Society},
  4(3):201--236, 2002.

\bibitem{ST02}
H~Mete Soner and Nizar Touzi.
\newblock Stochastic target problems, dynamic programming, and viscosity
  solutions.
\newblock {\em SIAM Journal on Control and Optimization}, 41(2):404--424, 2002.

\end{thebibliography}
\bibliographystyle{plain}
\end{document}